\documentclass[12pt]{article}
\usepackage[utf8]{inputenc}
\usepackage[T2A]{fontenc}
\usepackage[english]{babel}
\usepackage{amsfonts,amsmath, amssymb}
\usepackage{xcolor}

\usepackage{enumerate}
\usepackage{epsf,epsfig,amsfonts,a4wide}
\usepackage{amsfonts}
\usepackage{amsmath,amssymb,amsthm}
\usepackage{url}

\usepackage{amssymb}
\usepackage{amsthm}
\usepackage{epsf,epsfig,amsfonts,graphicx}
\usepackage{a4wide}
\newtheorem{Theorem}{Theorem}%[section]

\theoremstyle{definition}

\newtheorem{Remark*}{Remark}

\newcommand{\Ker}{\operatorname{Ker}}
\def\RR{\mathbb{R}}

\newcommand{\FF}{\mathbb{F}}

	 \newcommand{\KK}    {\mathbb K}
		\newcommand{\BB}    {\mathcal M}
   \newcommand{\cF}     {\mathcal F}
\newcommand{\N}{{\mathcal C}}
\newcommand{\cD}{{\mathcal D}}

\newcommand{\cB}{{\mathcal B}}

\def\cZ{\mathcal{Z}}
\def\cP{\mathcal{P}}
\newcommand{\extr}{{{extr}}}

\newcommand{\num}{{{num}}}
   \DeclareMathOperator {\dd} {d}
   \DeclareMathOperator {\Forg} {Forg}
   \DeclareMathOperator {\Diff} {Diff}

   \DeclareMathOperator {\Fix} {Fix}

   \renewcommand {\Im} {Im}

\title{Topology of spaces \\ of smooth functions and gradient-like flows \\ with prescribed singularities on surfaces}

\author{Elena A.\ Kudryavtseva
\thanks{Moscow State University, Moscow, Russia;
Moscow Center for Fundamental and Applied Mathematics, Moscow, Russia.
E-mail address: eakudr@mech.math.msu.su}
}
\date{}

\begin{document}
\maketitle

\begin{abstract}
By a gradient-like flow on a closed orientable surface $M$, we mean a closed 1-form $\beta$ defined on $M$ punctured at a finite set of points (sources and sinks of $\beta$) such that there exists a Morse function $f$ on $M$, called an energy function of $\beta$,
whose critical points coincide with equilibria of $\beta$, and the pair $(f,\beta)$ has a canonical form near each critical point of $f$. 
Let $\cB=\cB(\beta_0)$ be the space of all gradient-like flows on $M$ having the same types of local singularities as a flow $\beta_0$, and $\cF=\cF(f_0)$ the space of all Morse functions on $M$ having the same types of local singularities as an energy function $f_0$ of $\beta_0$. We prove that the spaces $\cF$ and $\cB$, equipped with $C^\infty$ topologies, are homotopy equivalent to some manifold $\BB_s$, moreover their decompositions into $\Diff^0(M)$-orbits are given by two transversal fibrations on $\BB_s$. Similar results are proved for topological equivalence classes on $\cF$ and $\cB$, and for non-Morse singularities.


{\bf MSC:} 58K05, 37E35, 37C15, 37B35, 37C86, 53C12, 58D27
%- 58K05   	Critical points of functions and mappings on manifolds
%- 37E35   	Flows on surfaces
%- 37C15   	Topological and differentiable equivalence, conjugacy, moduli, classification of dynamical systems
%- 37B35   	Gradient-like and recurrent behavior; isolated (locally maximal) invariant sets; attractors, repellers for topological dynamical systems
%- 37C86   	Foliations generated by dynamical systems
%- 53C12   	Foliations (differential geometric aspects) [See also 57R30, 57R32]
%- 58D27   	Moduli problems for differential geometric structures
\end{abstract}

%===============================================================================
%   Text of the theses
%===============================================================================

\section{Spaces of smooth functions with prescribed local singularities on surfaces} \label {sec:func}

Let $M$ be a smooth orientable connected closed two-dimensional surface, and $f_0\in C^\infty(M)$ a function all whose critical points have types $A_\mu$, $D_\mu$, $E_\mu$ (e.g.\ a Morse function). 

Recall that a point $P\in M$ is {\em critical} for $f\in C^\infty$ if $\dd f(P)=0$. A function $f\in C^\infty(M)$ is called {\em Morse} if all its critical points are {\em non-degenerate} (of type $A_1$), i.e.\ $\dd^2 f(P)$ is non-degenerate. By the Morse lemma, locally $f=\pm x^2\pm y^2+f(P)$ in suitable coordinates near each critical point $P$.

Consider the set $\mathcal{F}=\mathcal{F}(f_0)$ of all functions $f\in C^\infty(M)$ having the same types of local singularities as $f_0$.
Denote by $\mathcal{D}^0(M)$ the identity path component of the group $\mathcal{D}(M)=\mathrm{Diff}^+(M)$ of orientation-preserving diffeomorphisms endowed with $C^\infty$ topology. 
The group $\mathcal{D}(\RR)\times\mathcal{D}(M)$ acts on the space $\mathcal{F}$ by ``left-right changes of coordinates''.

We want to describe the topology of the space $\mathcal{F}$, equipped with the $C^\infty$-topology, and its decomposition into $\cD^0(\RR)\times\cD^0(M)$- and $\mathcal{D}^0(M)$-orbits. This problem was solved by the author in the cases when either $f_0$ is a Morse function and $\chi(M)<0$ \cite{K12,K12v}, or all critical points of $f_0$ have $A_\mu$ types, $\mu\in\mathbb{N}$ \cite{K16}. Topology of the $\mathcal{D}^0(M)$-orbits was studied by S.I.~Maksymenko \cite{Maks} (allowing some other types of degenerate singularities) and by the author \cite{K12,K12v,K16} (for $A_\mu$-singularities).

For any function $f\in\mathcal{F}$, consider the set $\mathcal{C}_f := \{P\in M\mid \mathrm{d} f(P)=0\}$ of its critical points. These critical points form five classes of topological equivalence (some classes may be empty):
$$
\mathcal{C}_f^{min} = \bigcup\limits_{i\ge1} A_{2i-1}^{+,+}(f), \
\mathcal{C}_f^{max} = \bigcup\limits_{i\ge1} A_{2i-1}^{+,-}(f), \
\mathcal{C}_f^{saddle} = A_1^-(f) \cup \bigcup\limits_{\eta=\pm}(\bigcup\limits_{i\ge2} A_{2i-1}^{-,\eta}(f) \cup D_{2i+1}^\eta(f))\cup E_7^\eta(f),
$$
$$
\mathcal{C}_f^{triv} = (\bigcup\limits_{i\ge1,\eta=\pm} A_{2i}^\eta(f))\cup (\bigcup\limits_{j\ge2} D_{2j}^+(f)) \cup E_6^+(f)\cup E_6^-(f) \cup E_8^+(f)\cup E_8^-(f), \quad 
\mathcal{C}_f^{mult} = \bigcup\limits_{j\ge2} D_{2j}^-(f),
$$
i.e.\ the critial points of local minima, local maxima, saddle points, quasi- and multy-saddle points, respectively. Here $A_\mu^{\pm,\pm}(f)$, $D_\mu^\pm(f)$ and $E_\mu^\pm(f)$ denote the corresponding subsets of critical points of $A-D-E$ types.
In the set $\mathcal{C}^{extr}_f:=\mathcal{C}^{min}_f\cup\mathcal{C}^{max}_f$ of local extremum points, consider the subset $\mathcal{C}^{{extr}*}_f$ of degenerate (non-Morse) critical points.

Denote $s:=\max\{0,\chi(M)+1\}$.

\begin{Theorem} \label {thm:1}
For any function $f_0\in C^\infty(M)$, whose all critical points have $A-D-E$ types (e.g.\ a Morse function), the space $\mathcal{F}=\mathcal{F}(f_0)$ has the homotopy type of a manifold $\BB_s=\BB_s(f_0)$ having dimension $\dim \BB_s=2s + |\mathcal{C}_{f_0}| + |\mathcal{C}^{{extr}*}_{f_0}| + |\mathcal{C}^{triv}_{f_0}| + 2|\mathcal{C}^{saddle}_{f_0}|+3|\mathcal{C}^{mult}_{f_0}|$. Moreover:
\begin{itemize}
\item[(a)] There exists a surjective submersion $\kappa:\mathcal{F}\to \BB_s$ and a stratification \cite{Whi65} (respectively, a fibration of codimension $|\mathcal{C}_{f_0}|$) on $\BB_s$ such that every $\mathcal{D}^0(\RR)\times\mathcal{D}^0(M)$-orbit (resp., $\mathcal{D}^0(M)$--orbit) in $\mathcal{F}$ is the $\kappa$-preimage of a stratum (resp., a fiber) in $\BB_s$.
\item[(b)] The map $\kappa$ provides a homotopy equivalence between every $\cD^0(M)$-invariant subset $I\subseteq\cF$ and its image $\kappa(I)\subseteq \BB_s$. In particular, it provides homotopy equivalences between $\cF$ and $\BB_s$, аnd between every $\cD^0(\RR)\times\cD^0(M)$--orbit (resp., $\cD^0(M)$--orbit) from $\cF$ and the corresponding stratum (resp., fiber) in $\BB_s$.
\end{itemize}
In particular, $\pi_k(\mathcal{F})\cong \pi_k(\BB_s),\ H_k(\mathcal{F})\cong H_k(\BB_s).$ Thus $H_k(\mathcal{F})=0$ for all $k>\dim \BB_s$.
\end{Theorem}

\begin{Remark*} \label {rem:1}
Denote by $\cF^1=\cF^1(f_0)$ the space of all functions $f\in\cF$ whose all local extrema equal $\pm 1$ and the sum of values at all non-extremum critical points vanishes. Then an analogue of Theorem \ref {thm:1} holds when $\cF$ and $\BB_s$ are replaced by $\cF^1$ and a submanifold $\BB^1_s\subset\BB_s$, respectively, where $\BB^1_s$ is a union of fibres of $\BB_s$, $\dim\BB^1_s=\dim\BB_s-|\N^\extr_{f_0}|-1
=2s + |\mathcal{C}^{{extr}*}_{f_0}| + 2|\mathcal{C}^{triv}_{f_0}| + 3|\mathcal{C}^{saddle}_{f_0}|+4|\mathcal{C}^{mult}_{f_0}| - 1$. Actually $\BB^1_s$ is a strong deformation retract of $\BB_s$, so it is homotopy equivalent to $\BB_s$.
\end{Remark*}

Our proof of Theorem \ref {thm:1} uses results \cite {K09, Ore21} about a ``uniform'' reduction of a smooth function to a normal form near its critical points.

\section{Morse flows and gradient-like flows on surfaces} \label {sec:flow}

Suppose $\Omega\in\Lambda^n(M)$ is a volume form on a $n$-manifold $M=M^n$. Let $\mathcal{P}\subset M$ be a finite subset. For any vector field $\xi$ on $M':=M\setminus\mathcal{P}$, we assign the $(n-1)$-form $\beta=i_\xi\Omega\in\Lambda^{n-1}(M')$. Clearly, this assignement is one-to-one, and $\xi\in\Ker\beta$. Furthermore, the flow of the vector field $\xi$ is volume-preserving if and only if $\beta$ is a closed form. Indeed: the Lie derivative $L_\xi\Omega=(i_\xi\mathrm{d}+\mathrm{d} i_\xi)\Omega=\mathrm{d} i_\xi\Omega=\mathrm{d}\beta$, so the Lie derivative vanishes if and only if $\mathrm{d}\beta=0$. By abusing language, we will call the $(n-1)$-form $\beta$ a {\em flow}. 

Suppose now that $n=\dim M=2$. Denote $\cZ_\beta:=\{P\in M'\mid \beta(P)=0\}$.

A closed 1-form $\beta$ on $M'=M\setminus\mathcal{P}$ will be called a {\em Morse flow} on $M$ if, in a neighbourhood of every point $P\in \mathcal{P}\cup\cZ_\beta$, there exist local coordinates $x,y$ such that either 
$\beta=\mathrm{d}(2xy)=\dd(\Im(z^2))$ and $P\in\mathcal{Z}_\beta$, or 
$\beta=\pm(x\mathrm{d} y-y\mathrm{d} x)/(x^2+y^2)=\pm\dd(\Im(\ln z))$ and $P\in\cP$, where $z=x+iy$.
Geometrically, the set $\mathcal{P}_\beta:=\mathcal{P}$ consists of {\em sources} and {\em sinks} of the flow $\beta$,
while the set $\mathcal{Z}_\beta$ consists of {\em saddle points} of the flow $\beta$.

A closed 1-form $\beta$ on $M'=M\setminus\mathcal{P}$ will be called a {\em gradient-like flow} on $M$ if there exists a Morse function $f\in C^\infty(M)$, called an {\em energy function} of $\beta$, such that 
\begin{itemize}
\item[(i)] the set $\mathcal{P}$ coincides with the set of local extremum points of $f$, 
\item[(ii)] the 2-form $\mathrm{d}f\wedge\beta|_{M\setminus\mathcal{C}_f}$ has no zeros and defines a positive orientation on $M$, 
\item[(iii)] in a neighbourhood of every point $P\in\mathcal{C}_f$, there exist local coordinates $x,y$ such that either \\
$f=f(P)+x^2-y^2$, $\beta=\mathrm{d}(2xy)$ and $P\in\mathcal{Z}=\mathcal{C}_f\setminus\mathcal{P}$, or \\
$f=f(P)\pm(x^2+y^2)$, $\beta=(x\mathrm{d} y-y\mathrm{d} x)/(x^2+y^2)$ and $P\in\mathcal{P}$.
\end{itemize}
Geometrically, the set $\mathcal{P}_\beta:=\mathcal{P}$ of {\em sources} and {\em sinks} of the flow $\beta$ coincides with the set of local extremum points of the energy function $f$, while the set $\mathcal{Z}_\beta:=\mathcal{Z}=\mathcal{C}_f\setminus\mathcal{P}$ of {\em saddle points} of the flow $\beta$ coincides with the set of saddle critical points of $f$.

Let $\beta_0$ be a Morse flow on $M$. Consider the set $\mathcal{B}=\mathcal{B}(\beta_0)$ of all gradient-like flows $\beta$ having the same types of local singularities as $\beta_0$ (in particular, $|\mathcal{Z}_\beta|=|\mathcal{Z}_{\beta_0}|$ and $|\mathcal{P}_\beta|=|\mathcal{P}_{\beta_0}|$).

First of all, let us characterize gradient-like flows among all 2D Morse flows. The following theorem is similar to a result by S.~Smale characterizing gradient-like flows among all Morse-Smale flows \cite {Sma61}.

\begin{Theorem} [Characterization of 2D gradient-like flows]  \label{thm:4}
Let $\beta_0$ be a Morse flow on $M$. Then:

\begin{itemize}
\item[(a)] The space $\cB=\cB(\beta_0)$ of gradient-like flows is nonempty if and only if $\beta_0$ has at least one sink and at least one source.

\item[(b)] A Morse flow $\beta$ is gradient-like if and only if 
\begin{itemize}
\item[(i)] $\beta$ has at least one sink and at least one source,
\item[(ii)] every separatrix of $\beta$ has two endpoints belonging to $\cZ_{\beta}\cup\cP_{\beta}$,
and 
\item[(iii)] there is no an oriented closed curve $P_1P_2\dots P_{k-1}P_k$ ($k\ge2$) formed by oriented separatrices of $\beta$, where 
$P_i\in\cZ_\beta$ and $P_k=P_1$,
\end{itemize}

\item[(c)] The space $\cB$ of gradient-like flows is open in the space of Morse flows and is $\cD(M)$-invariant. The orbit space 
$\cB_\num/\cD^0(M)$ is a $2|\cZ_{\beta_0}|$-dimensional manifold, where $\cB_\num$ is the space of gradient-like flows with enumerated sinks and sources.
\end{itemize}
\end{Theorem}

\section{Spaces of gradient-like flows on surfaces} \label {sec:grad}

We want to describe the topology of the space $\cB=\cB(\beta_0)$, equipped with the $C^\infty$-topology, and its decomposition into $\mathcal{D}^0(M)$-orbits and into classes of (orbital) topological equivalence.

\begin{Theorem} \label{thm:2}
For any gradient-like flow $\beta_0$ on $M$, the space $\cB=\cB(\beta_0)$ has the homotopy type of the manifold $\BB^1_s=\BB^1_s(f_0)$ from Theorem \ref{thm:1} and Remark 1, where $f_0$ is an energy function of $\beta_0$. Moreover:
\begin{itemize}
\item[(a)] There exists a surjective submersion $\lambda:\cB\to \BB^1_s$, a stratification and a 
$(|\cZ_{\beta_0}|+2s-1)$-dimensional fibration on $\BB^1_s$ such that every class of orbital topological equivalence (resp., $\cD^0(M)$-orbit) in $\cB$ is the $\lambda$-preimage of a stratum (resp., a fibre) from $\BB^1_s$.
\item[(b)] The map $\lambda$ provides a homotopy equivalence between every $\cD^0(M)$-invariant subset $I\subseteq\cB$ and its image $\lambda(I)\subseteq \BB^1_s$. In particular, it provides a homotopy equivalence between $\cB$ and $\BB^1_s$, as well as between every class of topological equivalence (resp., $\cD^0(M)$--orbit) in $\mathcal{B}$ and the corresponding stratum (resp., fibre) in $\BB^1_s$. 
\item[(c)] All fibres and strata in $\BB^1_s$ (and, thus, all classes of topological equivalence and all $\cD^0(M)$--orbits in $\cB$) are homotopy equivalent either to a point, or to $T^2$, or to $SO(3)/G$ or to $S^2$, in dependence on whether $\chi(M)<0$, or $\chi(M)=0$, or $\chi(M)\cdot|\mathcal{Z}_{\beta_0}|>0$, or $\chi(M)>0$ and $|\mathcal{Z}_{\beta_0}|=0$, respectively, where $G$ is a finite subgroup of $SO(3)$.
\end{itemize}
In particuar, $\pi_k(\cB)\cong \pi_k(\BB^1_s),\ H_k(\cB)\cong H_k(\BB^1_s)$. Thus $H_k(\cB)=0$ for all $k>\dim \BB^1_s$.
\end{Theorem}

\begin{Remark*}
The fibrations on the manifold $\BB^1_s$ in Theorem \ref {thm:1} and Theorem \ref {thm:2} are transversal to each other, and intersections of their fibres are $2s$-dimensional submanifolds diffeomorphic to the space $M^s\setminus\Delta$ of $s$-point configurations on $M$, where $\Delta=\cup\Delta_{ij}$, $\Delta_{ij} = \{(P_1,\dots,P_s)\in M^s\mid P_i=P_j\}$. 
Consider the topological space obtained from $\BB^1_s$ by contracting each such a $2s$-dimensional submanifold to a point. This space 
is known as the universal moduli space of real-normalized meromorphic 1-forms on $M$ \cite{Gru:Kri10}.
\end{Remark*}

\section{Describing the classifying manifolds $\BB_s$ and $\BB_s^1$} \label {sec:class}

The manifold $\BB_s=\BB_s(\beta_0)$ from Theorem \ref {thm:1} can be constructed as follows.

Let us consider the topological spaces
$$
\FF:=\{(f,\beta)\in\cF\times\cB \mid f\mbox{ is an energy function of }\beta\}, \qquad 
\FF^1:=\{(f,\beta)\in\FF \mid f\in\cF^1\}
$$
endowed with $C^\infty$ topology \cite{KP}.

Let us fix a $s$-point subset $N_s\subset M$, $|N_s|=s$.
Consider the subgroup $\cD_s(M):=\{\phi\in\cD(M)\mid N_s\subseteq\Fix(\phi)\}$ of the group $\cD(M)$ endowed with $C^\infty$ topology.
Denote by $\cD_s^0(M)$ the identity path component of the group $\cD_s(M)$.
Define the moduli spaces
$$
\BB_s=\BB_s(f_0):=\FF/\cD_s^0(M), \qquad
\BB_s^1=\BB_s^1(f_0):=\FF^1/\cD_s^0(M)
$$ 
endowed with quotient topology. 

One can show that $\BB_s^1$ is a strong deformation retract of $\BB_s$ \cite{KP}. Furthermore, by using a ``uniform'' reduction of a smooth function to a normal form near its critical points \cite{K09, Ore21}, one can prove that the forgetful maps 
\begin{equation}\label{eq:Forg}
\Forg_1:\FF\to\cF, \qquad 
\Forg_1|_{\FF^1}:\FF^1\to\cF^1, \qquad 
\Forg_2:\FF^1\to\cB
\end{equation}
are homotopy equivalences (cf.\ \cite{KP} for $\cF$ and Morse singularities, for other singularity types the proof is similar).

One can show that the group $\cD_s^0(M)$ acts freely on $\FF$. 
Since this group is contractible, we have homeomorphisms 
$$
\FF\approx\cD_s^0(M)\times\BB_s , \qquad 
\FF^1\approx\cD_s^0(M)\times\BB_s^1 .
$$
Therefore, the projections 
\begin{equation} \label {eq:moduli}
\FF \to \BB_s, \qquad 
\FF^1 \to \BB_s^1
\end{equation}
are homotopy equivalences.
It is easy to show that $\BB_s$ is a smooth manifold equipped with two transversal fibrations (whose fibres will be called ``horizontal'' and ``vertical'', respectively), and $\BB_s^1$ is its submanifold consisting of horizontal fibres.
Namely, each horizontal fibre is the $\Forg_1$-preimage of a $\cD^0(M)$-orbit in $\cF$, while 
each vertical fibre is the $\Forg_2$-preimage of a $\cD^0(M)$-orbit in $\cB$.

Now, for proving Theorem \ref {thm:1} and Remark \ref {rem:1}, one should consider the manifolds $\BB_s$ and $\BB_s^1$ fibred by the horizontal fibres.
For proving Theorem \ref {thm:2}, one should consider the manifold $\BB_s^1$ fibred by the vertical fibres. Then the theorems \ref {thm:1} and \ref {thm:2} follow from the homotopy equivalences \eqref{eq:Forg} and \eqref{eq:moduli}.
 
Each horizontal fibre (from Theorem \ref {thm:1} and Remark \ref{rem:1}) and each vertical fibre (from Theorem \ref {thm:2}) on the manifold $\BB^1_s$ are transversal to each other, and their intersection is a $2s$-dimensional submanifold diffeomorphic to $M^s\setminus\Delta$, the $s$-point configuration space on $M$, where $\Delta=\cup\Delta_{ij}$, $\Delta_{ij} = \{(P_1,\dots,P_s)\in M^s\mid P_i=P_j\}$.

\begin{Remark*}
Let us consider the moduli spaces
$$
\BB=\BB(f_0):=\FF/\cD^0(M), \qquad
\BB^1=\BB^1(f_0):=\FF^1/\cD^0(M)
$$ 
endowed with quotient topology. 
The space $\BB^1$ is known as the universal moduli space of real-normalized meromorphic differentials (or meromorphic 1-forms) on $M$ \cite{Gru:Kri10}.

Clearly, if $\chi(M)<0$ then $s=0$, $\BB_s=\BB_0=\BB$ and $\BB_s^1=\BB_0^1=\BB^1$, so $\BB$ and $\BB^1$ are manifolds.
If $\chi(M)\ge0$, then $\BB$ and $\BB^1$ are orbifolds in general, which can be obtained from the manifolds $\BB_s$ and $\BB_s^1$ (resp.) by contracting the intersection of each horizontal fibre with each vertical fibre to a point. 

Suppose that at least $s$ critical points of the function $f_0$ are enumerated (e.g., $\chi(M)<0$ or the singularity type of each of these $s$ points is different from the singularity type of any other critical point of $f_0$). Then:
\begin{itemize}
\item The orbifolds $\BB$ and $\BB^1$ are in fact manifolds, and there exist homeomorphisms
$$
\BB_s\approx (M^s\setminus\Delta)\times\BB, \qquad 
\BB_s^1\approx (M^s\setminus\Delta)\times\BB^1
$$ 
(cf.\ \cite{K12v, K13}). We remark that $M^s\setminus\Delta$ is homotopy equivalent to $\cD^0(M)$, which has the homotopy type of either a point or $T^2$ or $SO(3)$ in dependence on whether $\chi(M)<0$ or $\chi(M)=0$ or $\chi(M)>0$.
\item Each of the manifolds $\BB$ and $\BB^1$ (more precisely, their strong deformation retracts $\KK$ and $\KK^1$) is a ``skew cylindric-polyhedral complex'', i.e.\ it can be represented as the union of ``skew cylindric handles'' glued to each other in a nice way \cite{K12v, K13}. The skew cylindric handles of the manifold $\BB$ (resp., $\BB^1$) are in one-to-one correspondence with the $\cD^0(\RR)\times\cD^0(M)$-orbits in the space $\cF$ (resp., $\cF^1$). The Morse index of each handle equals the codimension of the corresponding $\cD^0(\RR)\times\cD^0(M)$-orbit in the space $\cF$ (resp., $\cF^1$).
\item Each skew cylindric handle of the manifold $\BB$ (resp., $\BB^1$) is incompressible, i.e.\ the inclusion mapping of the handle into the manifold $\BB$ (resp., $\BB^1$) induces a monomorphism of the fundamental groups.
\end{itemize}
\end{Remark*}

This work was supported by the Russian Foundation for Basic Research, grant No.~19-01-00775-a (the results of \S \ref {sec:func}) and the Russian Science Foundation, grant No.~17-11-01303 (the results of \S\S \ref {sec:flow}--\ref {sec:grad}).

%===============================================================================
%   Bibliography
%===============================================================================

\end{document}